\documentclass[12pt, reqno]{amsart}
\usepackage{amssymb}
\usepackage{amsfonts}
\usepackage{amssymb,amsmath,amscd}
\usepackage[colorlinks=true]{hyperref}
\usepackage{mathrsfs}
\newtheorem{Theorem}{\indent Theorem}[section]

\newtheorem{Lemma}[Theorem]{\indent Lemma}

\theoremstyle{remark}
\newtheorem{Remark}{Remark}

\begin{document}
\centerline{\bf On   a variant of the prime number theorem
}
\bigskip
\centerline{\small Wei Zhang}
\bigskip

\textbf{Abstract}
In this paper,  we can show that
 \begin{align*}
 S_{\Lambda}(x)=\sum_{1\leq n\leq x}\Lambda
 \left(\left[\frac{x}{n}\right]\right)=
\sum_{n=1}^{\infty}
\frac{\Lambda(n)}{n(n+1)}x
+O\left(x^{7/15+1/195+\varepsilon}\right),
 \end{align*}
where $\Lambda(n)$ is the von Mangdolt function. Moreover, we can also give   similar results related to the divisor function, which improve  previous results.

\textbf{Keywords}\ prime number theorem, exponential sum, divisor function

\textbf{2000 Mathematics Subject Classification}\  11N37, 11L07

\bigskip
\bigskip
\numberwithin{equation}{section}

\section{Introduction}
The prime number theory plays a central role in number theory, which sates that  there is a constant $c>0$ such that for $x\rightarrow\infty$, we have
\[
\sum_{1\leq n\leq x}\Lambda(n)=
x+O\left(x\exp(-c(\log x)^{3/5} (\log \log x)^{-1/5})\right)
\]
and the Riemann Hypothesis implies that
\[
\sum_{1\leq n\leq x}\Lambda(n)=
x+O\left(x^{1/2+\varepsilon}\right)
\]
where $\varepsilon$ is an arbitrarily small positive constant. Clearly it is also interesting to study
the distribution of prime numbers in different sequences of integers such as the arithmetic
progressions, the Beatty sequence $([\alpha n+\beta])_{n=1}^{\infty},$ the Piatetski-Shapiro sequence $([n^{c}])_{n=1}^{\infty}$,
etc, where $[t]$ denotes the integral part of the real number. For example, Banks and
Shparlinski (see Corollary 5.6 in \cite{BS}) proved the following result: Let $\alpha$ and $\beta$ be fixed real
numbers with $\alpha>1$, irrational and of finite type. Then there is a positive constant
$c=c(\alpha,\beta)$ such that
\[
\sum_{1\leq n\leq x}\Lambda([\alpha n+\beta])=x+O\left(x\exp(-c(\log x)^{3/5} (\log \log x)^{-1/5})\right).
\]
Piatetski-Shapiro sequences are named in honor of Piatetski-Shapiro, who proved that for any number $c\in(1,12/11)$ there are infinitely many primes of the form $[n^{c}]$ by showing that
 \begin{align}\label{01}
 \sum_{\substack{1\leq n\leq x\\ [n^{c}]\ \textup{is}\ \textup{prime}}}1
 =(1+o(1))\frac{N}{c\log N}.
 \end{align}
 The admissible range for $c$ in this problem has  been extended by many experts over the years.
And to date, the largest admissible $c$-range for (\ref{01}) seems to be $c\in(1,2817/2425)$ due to Rivat and Sargos  \cite{RS0} (see also the
references to the previous record holders they gave in their paper). Naturally, also
lower bound sieves have been employed, and the corresponding current record is a version of (\ref{01}) with a lower bound of the right order of magnitude instead of an
asymptotic formula and $c\in(1,243/205)$ due to Rivat and Wu \cite{RW}, which keeps for more than twenty years.
 We call the sum
\[
\sum_{1\leq n\leq x}f\left(\left[\frac{x}{n}\right]
\right)
\]
the fractional sum of $f$ (see \cite{St}), which has been considered in \cite{BDHPS}.
Similar to the well-known Beatty sequence $([\alpha n+\beta])_{n=1}^{\infty}$ and the Piatetski Shapiro sequence $([n^{c}])_{n=1}^{\infty}$, we may call the sequence
\[
\left(\left[\frac{x}{n}\right]
\right)_{n=1}^{\infty}
\]
a fractional sequence.
For this fractional sequence, Bordell\`{e}s-Dai-Heyman-Pan-Shparlinski \cite{BDHPS} established an asymptotic formula of
\[\sum_{1\leq n\leq x}f\left(\left[\frac{x}{n}\right]\right)
\]
under some simple assumptions of $f.$ Subsequently, Wu \cite{W01} and Zhai \cite{Zhai} improved their
results independently. In particular, applying  (Theorem 1.2(i) in \cite{W01}) or (Theorem 1 in \cite{Zhai}) to
the von Mangoldt function $\Lambda(n),$ we have
  \begin{align*}
 S_{\Lambda}(x)=\sum_{1\leq n\leq x}\Lambda
 \left(\left[\frac{x}{n}\right]\right)=
\sum_{n=1}^{\infty}
\frac{\Lambda(n)}{n(n+1)}x
+O\left(x^{1/2+\varepsilon}\right),
 \end{align*}
 where $\varepsilon$ is an arbitrary small positive constant.
With the help of the Vaughan's  identity and the method of one-dimensional
exponential sum, Ma and Wu \cite{MW} breaked the $1/2$-barrier:
  \begin{align*}
 S_{\Lambda}(x)=\sum_{1\leq n\leq x}\Lambda
 \left(\left[\frac{x}{n}\right]\right)=
\sum_{n=1}^{\infty}
\frac{\Lambda(n)}{n(n+1)}x
+O\left(x^{35/71+\varepsilon}\right).
 \end{align*}
Very recently, Bordell\`{e}s (see Corollary 1.3 in \cite{Bo}) sharpened the exponent $35/71$ to $97/203$ by using a
result of Baker \cite{Baker} on 2-dimensional exponential sums (see  Theorem 6 in \cite{Baker}).
Subsequently, by establishing an estimate
on 3-dimensional exponential sums, Liu-Wu-Yang \cite{LWY} proposed a better exponent $9/19$ (or $7774/16415$). The estimate
on 3-dimensional exponential sums of Liu-Wu-Yang \cite{LWY1} is a little different from the classical form.
In order to obtain much better result than the result of Liu-Wu-Yang, one need to generalize Robert and Sargos's results \cite{RS} to a special form, which is as done in Liu-Wu-Yang \cite{LWY1}, who proved the special 3-dimensional exponential sum by adapting and refining  Heath-Brown's approach (see \cite{HB}). A generalization of Robert and Sargos's results needs to adapting the approach of Fourvy-Iwaniec \cite{FI}. It seems that it is not easy to generalize   Robert and Sargos's results to our ideal form.
Hence for this  particular problem, the best possible exponent may be $9/19$ (or $7774/16415$) under the current technique.
The aim of this paper is to break  the  $9/19$-barrier (or $7774/16415$-barrier) and to further improve the results in \cite{Bo,MW,LWY1} by showing the following theorem.

\begin{Theorem}\label{th2}
Let $\Lambda(n)$ be the von Mangdolt function.
Then we have
 \begin{align*}
 S_{\Lambda}(x)=\sum_{1\leq n\leq x}\Lambda
 \left(\left[\frac{x}{n}\right]\right)=
\sum_{n=1}^{\infty}
\frac{\Lambda(n)}{n(n+1)}x
+O\left(x^{7/15+1/195+\varepsilon}\right),
 \end{align*}
 where $\varepsilon$ is an arbitrary small positive constant.
\end{Theorem}
\begin{Remark}
Note that $7/15+1/195\approx 0.47179$ and $9/19\approx 0.47368.$ If we choose $(\kappa,\lambda)=(1653/3494+\varepsilon,
1760/3494+\varepsilon)
=BA^{5}(13/84+\varepsilon,
55/84+\varepsilon)$ in applying Lemma \ref{LWY}, and by optimizing parameters,
one may obtain a little better error term. This basic observation can   be seen in \cite{LWY1}. Moreover, our results are obviously not the best by involving our ideas. Because we do not have optimization the parameters $N$ and $D$. Taking parameters optimally and considering the contribution of $H,$ one can yield better results. For example, in (\ref{remark}), in fact, our method implies an error term such that
\[
O\left(x^{7/15+32\varrho/45+\varepsilon}
+
x^{11/24+7\varpi/12+\varepsilon}+
x^{1/2-\varpi-\varrho+\varepsilon}
\right).
\]
By choosing $\varrho=6/923,$   $\varpi=20.5/923$ and $N=x^{7/15+64/13845},$ we can obtain such that
 \begin{align*}
 S_{\Lambda}(x)=\sum_{1\leq n\leq x}\Lambda
 \left(\left[\frac{x}{n}\right]\right)=
\sum_{n=1}^{\infty}
\frac{\Lambda(n)}{n(n+1)}x
+O\left(x^{7/15+64/13845+\varepsilon}\right).
 \end{align*}
 Let $(\kappa_{1},\lambda_{1})$ and $(\kappa_{2},\lambda_{2})$ be exponent pairs.
If we replace $3/8$ by $(3-3\lambda_{1})
/(5-3\lambda_{1}+\kappa_{1})$ in using the dichotomy and use the general form of Proposition 4.1 in \cite{LWY}, then we can obtain that for

\[
\varrho=\frac{(3/2-3\lambda_{1}/2)(8\kappa_{2}+2\lambda_{2}+14)
-(1-\lambda_{2})
(27-15\lambda_{1}+6\kappa_{1})}
{(47-27\lambda_{1}+10\kappa_{1})
(8\kappa_{2}+2\lambda_{2}+14)
-(27-15\lambda_{1}+6\kappa_{1})
(8\kappa_{2}+8)},
\]
\[
\varpi=\frac{(1-\lambda_{2})(47-27\lambda_{1}+10\kappa_{1})
-
(3/2-3\lambda_{1}/2)(8\kappa_{2}+8)
}{(47-27\lambda_{1}+10\kappa_{1})
(8\kappa_{2}+2\lambda_{1}+14)
-(27-15\lambda_{1}+6\kappa_{1})
(8\kappa_{2}+8)},
\]
and $N^{1/2-\varpi-\varrho},$ we have
 \begin{align*}
 S_{\Lambda}(x)=\sum_{1\leq n\leq x}\Lambda
 \left(\left[\frac{x}{n}\right]\right)=
\sum_{n=1}^{\infty}
\frac{\Lambda(n)}{n(n+1)}x
+O\left(x^{1/2-\varpi-\varrho+\varepsilon}\right).
\end{align*}
If we choose $(\kappa_{1},\lambda_{1})
=(\kappa_{2},\lambda_{2})=(1/2,1/2),$ then we can obtain that an error term $O(x^{7/15+64/13845}).$
If we choose $(\kappa_{1},\lambda_{1})
=(\kappa_{2},\lambda_{2})=(1653/3494+\varepsilon,
1760/3494+\varepsilon),$ then we can obtain   an error term $O(x^{0.471183811}).$
It is worth emphasizing that only by involving the ideas in this paper, the best possible result should be 7/15.
\end{Remark}
On the other hand, some experts also considered such type problems related to the divisor functions \cite{MS,LWY,Bo,St}.
By using the symmetry of the divisor function, in  \cite{MS},  it is proved that
\begin{align*}
 S_{\tau }(x)=\sum_{1\leq n\leq x}\tau \left(\left[\frac{x}{n}\right]\right)=
\sum_{n=1}^{\infty}
\frac{\tau (n)}{n(n+1)}x
+O\left(x^{11/23+\varepsilon}\right),
 \end{align*}
where $\tau (n)$ is  the number of representations of $n$ as product of two natural numbers and
\[
11/23\approx0.4782.
\]
Recently, this result was generalized and improved by many experts.
For example, Bordell\`es \cite{Bo} showed that
 \begin{align*}
 S_{\tau_{k}}(x)=\sum_{1\leq n\leq x}\tau_{k}\left(\left[\frac{x}{n}\right]\right)=
\sum_{n=1}^{\infty}
\frac{\tau_{k}(n)}{n(n+1)}x
+O\left(x^{\theta(k)+\varepsilon}\right),
 \end{align*}
 where $\theta(2)=19/40,$ $\theta(3)=283/574$ and
 \[
 \theta(k)=\frac12-\frac{1}{2(4k^{3}-k-1)}
 \]
 for $k\geq4$, where $\tau_{k}(n)$ is the number of representations of $n$ as product of $k$ natural numbers.
 Bordell\`es's result was improved by many experts. By  using a deep result of Jutila \cite{Ju},
 Stucky \cite{St} showed that $\theta(2)=5/11$ and Liu-Wu-Yang \cite{LWY} gave that
 \[
 \theta(k)=\frac{5k-1}{10k-1}
 \]
 for $k\geq 3.$
 Recently, in \cite{ZW}, the author showed that
 $\theta(k)=9/19$ for $k\geq 3.$ By   applying   similar arguments as the proof of Theorem \ref{th2} and the key ideas in \cite{ZW}, we can give the following improved result. Analogously, we also remark that we can obtain   a better error term $O(x^{0.471183811})$ (which is from Remark 1) for the following result.
\begin{Theorem}\label{th1}
Let $\tau_{k}(n)$ be the number of representations of $n$ as product of $k$ natural numbers.
Then we have
 \begin{align*}
 S_{\tau_{k}}(x)=\sum_{1\leq n\leq x}\tau_{k}\left(\left[\frac{x}{n}\right]\right)=
\sum_{n=1}^{\infty}
\frac{\tau_{k}(n)}{n(n+1)}x
+O\left(x^{7/15+1/195+\varepsilon}\right),
 \end{align*}
 where $\varepsilon$ is an arbitrary small positive constant.
\end{Theorem}
Our next aim of this paper is to consider the square-free divisor function over the fractional sequences.
The square-free divisor function has relation to the series
\[
\zeta^{2}(s)\zeta^{-1}(2s)
=\sum_{n=1}^{\infty}2^{\omega(n)}n^{-s},\ \ \Re(s)>1,
\]
where $\omega(n)$ denotes the number of distinct prime factors of $n$ and $\zeta(s)$ is the well known Riemann zeta function.
Let
$$\tau_{(2)}(n)=\sum_{\substack{d|n\\
d\in\mathcal{Q}_{2}}}1$$
 be the square-free divisor function, where $\mathcal{Q}_{2}$ is
the set of positive square-free integers and a number $q$ is called
square-free integer  if and only if
$
m^{2}|q\Longrightarrow m=1.
$
Then we have
$$\tau_{(2)}(n)=2^{\omega(n)}$$
and
\[
\tau_{(2)}(n)
=\sum_{n=n_{1}n_{2}}\tau(n_{1})
g(n_{2}),
\]
where $\tau(n)$ denotes the number of distinct factors of $n$ and
 \begin{align*}
g(d)=
 \begin{cases}
\mu(l)\ \ &\textup{if}\ d=l^{2},\\
\\
0\ \ &\textup{otherwise}.
\end{cases}
 \end{align*}

In \cite{Bo}, by using the Dirichlet hyperbolic method and more effort, it is proved that
\begin{align*}
 S_{\tau_{(2)}}(x)
=\sum_{1\leq n\leq x}\tau_{(2)}\left(\left[\frac{x}{n}\right]
\right)=
\sum_{n=1}^{\infty}
\frac{\tau_{(2)}(n)}{n(n+1)}x
+O\left(x^{97/202+\varepsilon}\right),
 \end{align*}
where  $[\cdot]$ denotes the floor function (i.e. the greatest integer function) and
$97/202\approx0.4802.$
Recently, this result was improved by many experts.
By using a new estimate on 3-dimensional exponential sums, in \cite{LWY}, Liu-Wu-Yang showed that
\begin{align*}
 S_{\tau_{(2)}}(x)
=\sum_{1\leq n\leq x}\tau_{(2)}\left(\left[\frac{x}{n}\right]
\right)=
\sum_{n=1}^{\infty}
\frac{\tau_{(2)}(n)}{n(n+1)}x
+O\left(x^{9/19+\varepsilon}\right),
 \end{align*}
 which has the same error term as their results for the primes over the fractional sequences.
We can improve $9/19\approx0.47368$ to $107/229\approx0.46724891$ by showing the following result.

\begin{Theorem}\label{th3}
Let $\tau_{(2)}(n)$ be the square-free divisor function.
Then we have
\begin{align*}
 S_{\tau_{(2)}}(x)
=\sum_{1\leq n\leq x}\tau_{(2)}\left(\left[\frac{x}{n}\right]
\right)=
\sum_{n=1}^{\infty}
\frac{\tau_{(2)}(n)}{n(n+1)}x
+O\left(x^{107/229+\varepsilon}\right),
 \end{align*}
 where $\varepsilon$ is an arbitrary small positive constant.
\end{Theorem}
\begin{Remark}
 Let $(\kappa_{1},\lambda_{1})$   be an exponent pair.
By using the dichotomy and   the general form of Lemma \ref{former11} in \cite{LWY1}, then we can obtain that for
\[
\varrho=\frac{38(1-\lambda_{1})
-(8+8\kappa_{1})}
{38(14+8\kappa_{1}+2\lambda_{1})
-22(8+8\kappa_{1})},
\]
\[
\varpi=\frac{ (14+8\kappa_{1}+2\lambda_{1})
-22(1-\lambda_{1}) }{38(14+8\kappa_{1}+2\lambda_{1})
-22(8+8\kappa_{1})},
\]
and $N^{1/2-\varpi-\varrho},$ we have
 \begin{align*}
 S_{\tau_{(2)}}(x)=\sum_{1\leq n\leq x}\tau_{(2)}
 \left(\left[\frac{x}{n}\right]\right)=
\sum_{n=1}^{\infty}
\frac{\tau_{(2)}(n)}{n(n+1)}x
+O\left(x^{1/2-\varpi-\varrho+\varepsilon}\right).
\end{align*}
If we choose $(\kappa_{1},\lambda_{1})
=(1/2,1/2),$ then we can obtain that an error term $O(x^{107/229+\varepsilon}).$
If we choose $(\kappa_{1},\lambda_{1})
=(1653/3494+\varepsilon,
1760/3494+\varepsilon),$ then we can obtain   an error term $O(x^{0.467135832}).$
\end{Remark}

\section{Preliminaries}

We will start  the proof for Theorem \ref{th2}  with some necessary lemmas.
The following lemma can be seen in Lemma 7 of Cao and Zhai \cite{CZ} or Theorem 2 in \cite{Baker1} (see also Lemma 3 of Baker-Harman-Rivat \cite {BHR}).
\begin{Lemma}\label{LWY}
Let $\alpha,$ $\beta,$ $\gamma$ be real numbers such that $\alpha(\alpha-1)\beta\gamma\neq 0.$
For $X>0,$ $H\geq1,$ $M\geq 1,$ and $N\geq 1,$ define
\[
S =S(H,M,N):=
\sum_{h\sim H}\sum_{m\sim M}\sum_{n\sim N}
a_{h,n}b_{m}e\left(\frac{X}{H^{\beta}
M^{\alpha}N^{\gamma}}
h^{\beta}m^{\alpha}n^{\gamma}\right),
\]
where $e(t)=e^{2\pi i t},$ the $a_{h,n}$ and $b_{m}$ are complex  numbers such that $a_{h,n}\leq 1,$ $b_{m}\leq 1$ and $m\sim M$ means that $M<m\leq 2M.$ For any $\varepsilon>0,$ we have
\[
S\ll
\left(\left(X^{\kappa}
H^{2+\kappa}M^{1+\kappa+\lambda}
N^{2+\kappa}\right)^{1/(2+2\kappa)}
+HM^{1/2}N+H^{1/2}MN^{1/2}+X^{-1/2}HMN\right)
X^{\varepsilon}
\]
uniformly for $M\geq 1,$ $N\geq 1,$   where
$(\kappa.\lambda)$ is an exponent pair and the implied constant may depend on $\alpha,\beta,\gamma,$ and $\varepsilon.$
\end{Lemma}
We need
the following  well-known lemma (for example, one can refer to page 441 of \cite {BB} or page 34 of  \cite{GK}).
\begin{Lemma}\label{z2}
Let $g^{(l)}(x)\asymp YX^{1-l}$ for $1< X\leq x\leq 2X$ and $l=1,2,\cdots.$ Then one has
\[
\sum_{X<n\leq 2X}e(g(n))\ll Y^{\kappa}X^{\lambda}+Y^{-1},
\]
where $(\kappa,\lambda)$ is any exponent pair.
\end{Lemma}

Let $\psi(t)=t-[t]-1/2$ for $t\in\mathbb{R}$ and $\delta\geq0.$
We also need
the following  well-known lemma.
This lemma can be seen in Theorem A.6 in \cite{GK} or Theorem 18 in \cite{Va}.
\begin{Lemma}\label{z3}
For $0<|t|<1,$ let
$$W(t) = \pi t(1-|t|)\cot\pi t + |t|.$$
  For $x\in\mathbb{R},$ $H\geq1,$ we define
$$\psi^{*}(x)=\sum_{1\leq |h|\leq H}(2\pi ih)^{-1}W\left(\frac{h}{H+1}\right)e(hx)$$
and
\[
\delta(x)=\frac{1}{2H+2}\sum_{|h|\leq H}\left(1-\frac{|h|}{H+1}\right)e(hx).
\]
Then $\delta(x)$ is non-negative, and we have
$$|\psi^{*}(x)-\psi(x)|\leq
\delta(x).$$
\end{Lemma}

To deal with the von Mangoldt function, we also need the following Vaughan's identity (for example, see \cite{LWY1} and the references therein).
\begin{Lemma}\label{pv} There are six real arithmetical functions $\alpha_{k}(n)$ verifying $\alpha_{k}(n)\ll _{\varepsilon} n^{\varepsilon}$
for
($n>1, 1\leq k\leq 6$) such that, for all $D>100$ and any arithmetical function $g,$ we have
\[\sum_{D<d\leq 2D}
\Lambda(d)g(d) = S_{1} +  S_{2} +  S_{3} +  S_{4},\]
where
\begin{align*}
&S_{1}=\sum_{m\leq D^{1/3}}\alpha_{1}(m)\sum_{D<mn\leq 2D}g(mn),\\
&S_{2}=\sum_{m\leq D^{1/3}}\alpha_{2}(m)\sum_{D<mn\leq 2D}g(mn)\log n,\\
&S_{3}=\mathop{\sum\sum}_{\substack{D^{1/3}<m,n\leq D^{2/3}\\ D<mn\leq 2D}}\alpha_{3}(m)\alpha_{4}(n)g(mn),
\\&S_{4}=\mathop{\sum\sum}_{\substack{D^{1/3}<m,n\leq D^{2/3}\\ D<mn\leq 2D}}\alpha_{5}(m)\alpha_{6}(n)g(mn).
\end{align*}
The sums $S_{1}$ and $S_2$ are called as type $I,$ $S_3$ and $S_4$ are called as type $II.$
\end{Lemma}
Next, we introduce a well-known result in \cite{RS}.
\begin{Lemma}\label{rs}
Let $\alpha,$ $\beta,$ $\gamma$ be real numbers such that $\alpha(\alpha-1)\beta\gamma\neq 0.$
For $X>0,$ $H\geq1,$ $M\geq 1,$ and $N\geq 1,$ define
\[
S =S(H,M,N):=
\sum_{h\sim H}\sum_{m\sim M}\sum_{n\sim N}
a_{h,n}b_{m}e\left(\frac{X}{H^{\beta}
M^{\alpha}N^{\gamma}}
h^{\beta}m^{\alpha}n^{\gamma}\right),
\]
where $e(t)=e^{2\pi i t},$ the $a_{h,n}$ and $b_{m}$ are complex  numbers such that $a_{h,n}\leq 1,$ $b_{m}\leq 1$ and $m\sim M$ means that $M<m\leq 2M.$ For any $\varepsilon>0,$ we have
\begin{align*}
S(H,M,N)(XHMN)^{-\varepsilon}
&\ll
\left(X
M^{2}N^{3}
H^{3}\right)^{1/4}
\\&+M\left(HN\right)^{3/4}
 +M^{1/2}HN
 +X^{-1/2}HNM,
\end{align*}
 where
the implied constant may depend on $\alpha,$ $\beta,$ $\gamma,$ and $\varepsilon.$
\end{Lemma}
The following lemma can be seen by (4.3) of Proposition 4.1 in \cite{LWY}, which is also needed in our proof.
\begin{Lemma}\label{former}
Let $\delta\notin-\mathbb{N}$ be a fixed constant. For $x^{6/13}\leq D \leq x^{2/3},$ we have
\[
\sum_{d\sim D}\Lambda(d)\psi\left(\frac{x}{d+\delta}\right)
\ll (x^{2}D^{7})^{1/12}x^{\varepsilon}.
\]
\end{Lemma}
 \begin{Lemma}\label{rs11}
For real numbers $\alpha_{1}, \alpha_{2}, \alpha_{3}$ such that $\alpha_{1}\alpha_{2}\alpha_{3}
(\alpha_{1}-1)(\alpha_{2}-2)\neq0$.
For $X>0,$ $M_{1}\geq1,$ $M_{2}\geq 1,$ and $M_{3}\geq 1,$ let
\[
S(M_{1},M_{2},M_{3}):=
\sum_{m_{2}\sim M_{2}}\sum_{m_{3}\sim M_{3}}
\left|\sum_{m_{1}\sim M_{1}}e\left(X\frac{m_{1}^{\alpha_{1}} m_{2}^{\alpha_{2}}
 m_{3}^{\alpha_{3}}}
{M_{1}^{\alpha_{1}}M_{2}^{\alpha_{2}}
M_{3}^{\alpha_{3}}}
\right)\right|,
\]
where $e(t)=e^{2\pi i t}.$   For any $\varepsilon>0,$ we have
\begin{align*}
S(M_{1},M_{2},M_{3})(XM_{1}M_{2}M_{3})^{-\varepsilon}&\ll
 \left(X
M_{1}^{2}M_{2}^{3}
M_{3}^{3}\right)^{1/4} +M_{1}^{1/2} M_{2}M_{3}
 +X^{-1}M_{1}M_{2}M_{3},
\end{align*}
 where
the implied constant may depend on $\alpha_{1},\alpha_{2},\alpha_{3},$ and $\varepsilon.$
\end{Lemma}
The following lemma can be seen  in \cite{LWY1}, which is proven by using the well the new 3-dimensional exponential sums of \cite{LWY} and also needed in our proof.
\begin{Lemma}\label{former11}
Let $\delta\notin-\mathbb{N}$ be a fixed constant. For $1\leq D \leq x^{8/11},$ we have
\[
\sum_{d\sim D}\tau_{(2)}(d)\psi\left(\frac{x}{d+\delta}\right)
\ll (x^{2}D^{7})^{1/12}x^{\varepsilon}.
\]
\end{Lemma}
\section{Proof of Theorem \ref{th2}}
Let
\[
\mathcal{N}=x^{7/15}
\]
We can write
\[
S_{\Lambda}(x):=S_{\Lambda,1} +S_{\Lambda,2} ,
\]
where
\begin{align}\label{c11}
S_{\Lambda,1}=\sum_{1\leq n\leq \mathcal{N}}\Lambda
\left(\left[\frac{x}{n}\right]\right)
\end{align}
and
\begin{align}\label{c12}
S_{\Lambda,2}=\sum_{\mathcal{N}<n\leq x}\Lambda
\left(\left[\frac{x}{n}\right]\right).
\end{align}
Obviously, by   $\Lambda(n)\ll n^{\varepsilon},$ we have
\begin{align*}
S_{\Lambda,1}=\sum_{1\leq n\leq \mathcal{N}}\Lambda\left(\left[\frac{x}{n}\right]\right)
&=\sum_{1\leq n\leq \mathcal{N}}(x/n)^{\varepsilon}
\\
&\ll \mathcal{N}^{1+\varepsilon}\\
&\ll x^{7/15+\varepsilon}.
\end{align*}

As to $S_{\Lambda,2},$  by    $\Lambda(n)\ll n^{\varepsilon},$
  we have
$$
\sum_{1\leq n\leq x} \Lambda(n)  \ll x^{1+\varepsilon}.
$$
Hence we can get
\begin{align}\label{c1}
\begin{split}
S_{\Lambda,2}&=\sum_{\mathcal{N}<n\leq x}\Lambda\left(\left[\frac{x}{n}\right]\right)
\\&
=\sum_{d\leq x/\mathcal{N}}\Lambda(d)
\sum_{x/(d+1)<n\leq x/d}1
\\&
=\sum_{d\leq x/\mathcal{N}}
\Lambda(d)\left(\frac{x}{d}-\frac{x}{d+1}
-\psi(\frac{x}{d})
+\psi(\frac{x}{d+1})\right)
\\&=x\sum_{d=1}^{\infty}
\frac{\Lambda(d)}{d(d+1)}
+O\left(
\mathcal{N}^{1+\varepsilon}
\right)\\
&+O\left((\log x)\max_{\mathcal{N}<D\leq x^{1/2+\varpi}}\left|\sum_{D<d\leq 2D}\Lambda(d)
\psi\left(\frac{x}
{d+\delta}\right)\right|\right)
\\
&+O\left((\log x)\max_{x^{ 1/2+\varpi}<D\leq x/\mathcal{N}}\left|\sum_{D<d\leq 2D}\Lambda(d)
\psi\left(\frac{x}{d+\delta}
\right)\right|\right),
\end{split}\end{align}
where   $\delta\in\{0,1\}.$
We need to consider $\mathcal{N}<D\leq x^{1/2+\varpi}$ and $x^{1/2+\varpi}<D\leq x/\mathcal{N}$ respectively. Here we restrict the range of $\varpi$ for $\varpi\in(0,1/38).$
If $\vartheta=1/38,$ we can obtain the error term $O(x^{9/19+\varepsilon})$ given by Liu-Wu-Yang \cite{LWY1}. We hope to balance the parameters by using $\vartheta<1/38$ to obtain a better result.
For $\mathcal{N}<D\leq x^{1/2+\varpi},$ we need to follow the arguments in \cite{LWY1}. By Lemma \ref{former}, we can  obtain that
\[
\sum_{D<d\leq 2D}\Lambda(d)
\psi\left(\frac{x}
{d+\delta}\right)
\ll x^{1/6}D^{7/12}\ll x^{11/24+7\varpi/12+\varepsilon}.
\]
Then for $x^{1/2+\varpi}<D\leq x/\mathcal{N}$, we need to estimate
\[\sum_{D<d\leq 2D}\Lambda(d)\psi
\left(\frac{x}{d+\delta}\right).
\]
By Lemma \ref{z3}, we have
\begin{align}\label{c2}
\begin{split}\sum_{D<d\leq 2D}&\Lambda(d)\psi
\left(\frac{x}{d+\delta}\right)\\
&\ll
\left|\sum_{1\leq h\leq H}\frac1h\sum_{D<d\leq 2D}\Lambda(d) e\left(\frac{hx}{d+\delta}\right)\right|
\\&+ \left|\sum_{1\leq h\leq H}\frac1H\sum_{D<d\leq 2D}\Lambda(d)
e\left(\frac{hx}{d+\delta}\right)\right|
+D/H.
\end{split}\end{align}
Then we will focus on the estimate of
\begin{align*}
\sum_{1\leq h\leq H}\frac1h\sum_{D<d\leq 2D}\Lambda(d)
e\left(\frac{hx}{d+\delta}\right).
\end{align*}
And we can  handle the sum
\[
\sum_{1\leq h\leq H}\frac1H\sum_{D<d\leq 2D}\Lambda(d)
e\left(\frac{hx}{d+\delta}\right)
\]
similarly.
By partial summation, we have
\begin{align}\label{n5}
\begin{split}
&\sum_{D<d\leq 2D}\sum_{1\leq h\leq H}\frac{\Lambda(d)}{h}
e\left(\frac{hx}{d}\right)
e\left(\frac{-\delta hx}{d(d+\delta)}\right).
\\&
\ll \int_{D}^{2D}e\left(\frac{-\delta hx}{u(u+\delta)}\right)d\left(
\sum_{D<d\leq u}\sum_{1\leq h\leq H}\frac{\Lambda(d)}{h}
e\left(\frac{hx}{d}\right)\right)\\
&\ll
\max_{D\leq D_{1}\leq 2D}\left|\sum_{D<d\leq D_{1}}\sum_{1\leq h\leq H}\frac{\Lambda(d)}{h}
e\left(\frac{hx}{d}\right)\right|\\
&+\frac{Hx}{D^{2}}\max_{D\leq D_{1}\leq 2D}\left|\sum_{D<d\leq D_{1}}\sum_{1\leq h\leq H}\frac{\Lambda(d)}{h}
e\left(\frac{hx}{d}\right)\right|.
\end{split}
\end{align}
Choose $H=D^{2}/x^{1-\varrho}.$ Here we restrict the range of $\varrho$ for $\varrho\in(0,1/4).$ Then one can verify that $H\geq 1.$
Hence, we only need to estimate the 
sum
 \begin{align*}
\max_{1\leq H_{1}\leq H}\sum_{h\sim H_{1}}\frac1h\sum_{D<d\leq D_{1}}\Lambda(d)
e\left(\frac{hx}{d}\right).
\end{align*}
Further, we need to use Lemma \ref{pv} to deal with this sum.
By Lemma \ref{pv}, there   are six real arithmetical functions $\beta_{k}(n)$ verifying $\beta_{k}(n)\ll _{\varepsilon} n^{\varepsilon}$
for
($n>1,\  1\leq k\leq 6$)  such that
\[
\sum_{h\sim H_{1}}\frac1h\sum_{D<d \leq D_{1}}\Lambda(d)e\left(\frac{hx}{d}\right)\ll
T_{1}+T_{2}+T_{3}+T_{4},
\]
where
\[
T_{1}\ll \sum_{h\sim  H_{1}
}\frac1h\sum_{m\leq D^{1/3}}\beta_{1}(m)\sum_{D<mn\leq D_{1}}e\left(\frac{hx}{mn}\right),
\]
\[
T_{2}\ll \sum_{h\sim  H_{1}
}\frac1h\sum_{m\leq D^{1/3}}\beta_{2}(m)\sum_{D<mn\leq D_{1}}e\left(\frac{hx}{mn}\right)\log n,
\]
\[
T_{3}\ll \sum_{h\sim  H_{1}
}\frac1h
\mathop{\sum\sum}_{\substack{D^{1/3}
<m,n\leq D^{2/3}\\ D<mn\leq D_{1}}}\beta_{3}(m)\beta_{4}(n)
e\left(\frac{hx}{mn}\right),
\]
and
\[
T_{4}\ll \sum_{h\sim  H_{1}
}\frac1h
\mathop{\sum\sum}_{\substack{D^{1/3}<m,n\leq D^{2/3}\\ D<mn\leq D_{1}}}\beta_{5}(m)\beta_{6}(n)
e\left(\frac{hx}{mn}\right).\]
For $T_{1}$, we have
\begin{align*}
T_{1}&\ll \sum_{h\sim  H_{1}
}\frac1h\sum_{1\leq m\leq D^{1/4}}\beta_{1}(m)\sum_{D<mn\leq D_{1}}e\left(\frac{hx}{mn}\right)\\
&+ \sum_{h\sim  H_{1}
}\frac1h\sum_{D^{1/4}\leq m\leq D^{1/3}}\beta_{1}(m)\sum_{D<mn\leq D_{1}}e\left(\frac{hx}{mn}\right)
\end{align*}
By Lemma \ref{z2} and choosing $(\kappa,\lambda)=(1/2,1/2)$, we have
\begin{align*}
&\sum_{h\sim  H_{1}
}\frac1h\sum_{1\leq m\leq D^{1/4}}\beta_{1}(m)\sum_{D<mn\leq D_{1}}e\left(\frac{hx}{mn}\right)\\
&\ll D^{\varepsilon}
\left(\frac{H^{1/2}x^{1/2}}{D^{1/4}}
+D^{2}/x\right)
\\&\ll x^{0.4+4\varrho/15+\varepsilon},
\end{align*}
where  we have used   $H=D^{2}/x^{1-\varrho}.$
By Lemma \ref{LWY}, we have
\begin{align*}
&\sum_{h\sim  H_{1}
}\frac1h\sum_{D^{1/4}\leq m\leq D^{1/3}}\beta_{1}(m)\sum_{D<mn\leq D_{1}}e\left(\frac{hx}{mn}\right)\\
&\ll D^{\varepsilon}
\left(x^{1/6}
D^{\frac23\times\frac{8}{15}\times\frac13}+ D^{7/8}+\frac{D^{3/2}}{x^{1/2}}\right)\\
&\ll x^{7/15+\varepsilon},
\end{align*}
Then   by choosing $\mathcal{N}=x^{7/15}$ and similar arguments for $T_{2},$ we have
\[
T_{1}+T_{2}  \ll x^{7/15+\varepsilon}.
\]

As for $T_{3}$ and $T_{4},$ by the symmetry, we have
\begin{align*}
T_{3}&\ll \sum_{h\sim  H_{1}
}\frac1h
\mathop{\sum\sum}_{
\substack{D^{1/3}<m\leq D^{1/2}\\
D^{1/2}<n\leq D^{2/3}\\
 D<mn\leq D_{1}}}\beta_{3}(m)\beta_{4}(n)
e\left(\frac{hx}{mn}\right)
 \\&+
\sum_{h\sim  H_{1}
}\frac1h
\mathop{\sum\sum}_{
\substack{D^{1/3}<n\leq D^{1/2}\\
D^{1/2}<m\leq D^{2/3}\\
 D<mn\leq D_{1}}}\beta_{3}(m)
 \beta_{4}(n)e\left(\frac{hx}{mn}\right)\\
 &:= U_{1}+U_{2},
\end{align*}
and
\begin{align*}
T_{4}&\ll
\sum_{h\sim  H_{1}
}\frac1h
\mathop{\sum\sum}_{
\substack{D^{1/3}<m\leq D^{1/2}\\
D^{1/2}<n\leq D^{2/3}\\
D<mn\leq D_{1}}}\beta_{5}(m)
 \beta_{6}(n)e\left(\frac{hx}{mn}\right)\\
&+
\sum_{h\sim  H_{1}
}\frac1h
\mathop{\sum\sum}_{
\substack{D^{1/3}<n\leq D^{1/2}\\
D^{1/2}<m\leq D^{2/3}\\
D<mn\leq D_{1}}}\beta_{5}(m)
 \beta_{6}(n)e\left(\frac{hx}{mn}\right)\\
 &:= U_{3}+U_{4}.
\end{align*}

Then we focus on the estimate of $U_{1}$ and other cases are similar by involving the symmetry.
We divide $U_{1}$ into two cases
\begin{align*}
U_{1}&\ll \sum_{h\sim  H_{1}
}\frac1h
\mathop{\sum\sum}_{
\substack{D^{1/3}<m\leq D^{3/8}\\
D^{5/8}<n\leq D^{2/3}\\
 D<mn\leq D_{1}}}\beta_{3}(m)\beta_{4}(n)
e\left(\frac{hx}{mn}\right)\\
 &+
 \sum_{h\sim  H_{1}
}\frac1h
\mathop{\sum\sum}_{
\substack{D^{3/8}<m\leq D^{1/2}\\
D^{1/2}<n\leq D^{5/8}\\
 D<mn\leq D_{1}}}\beta_{3}(m)\beta_{4}(n)
 e\left(\frac{hx}{mn}\right)\\
 & := U_{11}+U_{12}.
\end{align*}
By Lemma \ref{LWY}, choosing $(\kappa,\lambda)=(1/2,1/2),$ we have
\begin{align*}
U_{11}&\ll
D^{\varepsilon}\sum_{h\sim  H_{1}
}\frac1h
\mathop{\sum\sum}_{
\substack{m\sim M\\
n\sim N\\
 D<mn\leq D_{1}}}\beta_{3}(m)\beta_{4}(n)
e\left(\frac{hx}{mn}\right)\\
&\ll D^{\varepsilon}\left(\left(x^{1/2}M^{2}N^{3/2}\right)^{1/3}
+M^{1/2}N+MN^{1/2}+\frac{M^{3/2}N^{3/2}}
{x^{1/2}}\right),
\end{align*}
where $D^{1/3}\leq M \leq D^{3/8}$ and $D^{5/8}\leq N\leq D^{2/3}.$
Recall that $\mathcal{N}=x^{7/15}$ and $D\leq x/\mathcal{N},$ we have
\begin{align*}
U_{11}&\ll
D^{\varepsilon}\sum_{h\sim  H_{1}
}\frac1h
\mathop{\sum\sum}_{
\substack{m\sim M\\
n\sim N\\
 D<mn\leq D_{1}}}\beta_{3}(m)\beta_{4}(n)
e\left(\frac{hx}{mn}\right)\\
&\ll
D^{\varepsilon}
\left((x^{1/2}D^{3/2+3/16})^{1/3}
+D^{5/6}+D^{3/2}/x^{1/2}\right)
\\
&\ll x^{7/15+\varepsilon}.
\end{align*}
By Lemma \ref{rs},  we have
\begin{align*}
U_{12}&\ll
D^{\varepsilon}\sum_{h\sim  H_{1}
}\frac1h
\mathop{\sum\sum}_{
\substack{m\sim M\\
n\sim N\\
 D<mn\leq D_{1}}}\beta_{3}(m)\beta_{4}(n)
e\left(\frac{hx}{mn}\right)\\
&\ll D^{\varepsilon}
\left(\left(xMN^{2}\right)^{1/4}
+MN^{3/4}
+M^{1/2}N+\frac{M^{3/2}N^{3/2}}
{x^{1/2}}\right),
\end{align*}
where $D^{3/8}\leq M \leq D^{1/2}$ and $D^{1/2}\leq N\leq D^{5/8}.$
Recall that $\mathcal{N}=x^{7/15}$ and $D\leq x/\mathcal{N},$ we have
\begin{align*}
U_{12}&\ll
D^{\varepsilon}\sum_{h\sim  H_{1}
}\frac1h
\mathop{\sum\sum}_{
\substack{m\sim M\\
n\sim N\\
 D<mn\leq D_{1}}}\beta_{3}(m)\beta_{4}(n)
e\left(\frac{hx}{mn}\right)\\
&\ll
D^{\varepsilon}
\left(\left(xD^{2-3/8}\right)^{1/4}
+D^{7/8}
+\frac{M^{3/2}N^{3/2}}
{x^{1/2}}\right)
\\
&\ll x^{7/15+\varepsilon}.
\end{align*}
Then we have
\[
U_{1}\ll U_{11}+U_{12}\ll x^{7/15+\varepsilon}.
\]
Similarly, we can obtain the estimates of $U_{2}, U_{3}, U_{4}.$ Hence we have
\[
T_{3}+T_{4}\ll x^{7/15+\varepsilon}.
\]
Recall the estimates of $T_{1}$ and $T_{2}$, we can obtain that
\begin{align*}
\sum_{h\sim H_{1}}\frac1h\left|\sum_{D<n\leq D_{1}}\Lambda(n)e\left(\frac{hx}{n}\right)\right|
 \ll x^{7/15+\varepsilon}.
\end{align*}
Then by (\ref{c1})-(\ref{n5}), we have
\begin{align}\label{remark}
\begin{split}
S_{\Lambda,2}&=\sum_{N<n\leq x}\Lambda
\left(\left[\frac{x}{n}\right]\right)
=x\sum_{d=1}^{\infty}\frac{\Lambda(d)}{d(d+1)}
\\
&+O\left(x^{7/15+\varrho+\varepsilon}
+
x^{11/24+7\varpi/12+\varepsilon}+
x^{1/2-\varpi-\varrho+\varepsilon}
\right).
\end{split}
\end{align}
Choose
$\varrho=1/195$ and $\varpi=3/130.$
Recall that
\[
S_{\Lambda,1}=\sum_{1\leq n\leq \mathcal{N}}\Lambda
\left(\left[\frac{x}{n}\right]\right)
\ll x^{7/15+\varepsilon}.
\]
Then by (\ref{c11}) and (\ref{c12}), we have
\[
\sum_{1<n\leq x}\Lambda
\left(\left[\frac{x}{n}\right]\right)
=x\sum_{d=1}^{\infty}\frac{\Lambda(d)}{d(d+1)}
+O(x^{7/15+1/195+\varepsilon}).
\]
This completes the proof.


\section{Proof of Theorem \ref{th1}}
Now we begin the proof of Theorem \ref{th1}. Let
\[
\mathcal{N}=x^{7/15}.
\]
We can write
\[
S_{\tau_{k}}(x):=S_{\tau_{k},1} +S_{\tau_{k},2} ,
\]
where
\begin{align}\label{c110}
S_{\tau_{k},1}=\sum_{n\leq \mathcal{N}}\tau_{k}
\left(\left[\frac{x}{n}\right]\right)
\end{align}
and
\begin{align}\label{c120}
S_{\tau_{k},2}=\sum_{\mathcal{N}<n\leq x}\tau_{k}
\left(\left[\frac{x}{n}\right]\right).
\end{align}
Obviously, by   $\tau_{k}(n)\ll n^{\varepsilon},$ we have
\begin{align*}
S_{\tau_{k},1}=\sum_{n\leq \mathcal{N}}
\tau_{k}\left(\left[\frac{x}{n}
\right]\right)
&=\sum_{n\leq \mathcal{N}}(x/n)^{\varepsilon}
\\
&\ll \mathcal{N}^{1+\varepsilon}\\&\ll x^{7/15+\varepsilon}.
\end{align*}
As to $S_{\tau_{k},2},$  by    $\tau_{k}(n)\ll n^{\varepsilon},$
  we have
$$
\sum_{n\leq x} \tau_{k}(n)  \ll x^{1+\varepsilon}.
$$
Hence we can get
\begin{align}\label{c100000}
\begin{split}
S_{\tau_{k},2}&=\sum_{\mathcal{N}<n\leq x}\tau_{k}\left(\left[\frac{x}{n}\right]\right)
\\&
=\sum_{d\leq x/\mathcal{N}}
\tau_{k}(d)\sum_{x/(d+1)<n\leq x/d}1
\\&
=\sum_{d\leq x/\mathcal{N}}\tau_{k}(d)
\left(\frac{x}{d}-\frac{x}{d+1}
-\psi(\frac{x}{d})
+\psi(\frac{x}{d+1})\right)
\\&=x\sum_{d=1}^{\infty}
\frac{\tau_{k}(d)}{d(d+1)}
+O\left(
\mathcal{N}^{1+\varepsilon}
\right)\\
&+O\left((\log x)\max_{\mathcal{N}<D\leq x^{1/2+\varpi}}\left|\sum_{D<d\leq 2D}\tau_{k}(d)
\psi\left(\frac{x}{d+\delta}\right)
\right|\right)
\\
&+O\left((\log x)\max_{x^{1/2+\varpi}<D\leq x/\mathcal{N}}\left|\sum_{D<d\leq 2D}\tau_{k}(d)
\psi\left(\frac{x}{d+\delta}\right)\right|\right),
\end{split}\end{align}
where $\mathcal{N}\leq D\leq x/\mathcal{N}$ and $\delta\in\{0,1\}.$ We need to consider $\mathcal{N}<D\leq x^{ 1/2+\varpi}$ and $x^{1/2+\varpi}<D\leq x/\mathcal{N}$ respectively.
For $\mathcal{N}<D\leq x^{1/2+\varpi},$ we need to follow the arguments in \cite{LWY,LWY1,ZW} to obtain that
\[
\sum_{D<d\leq 2D}\tau_{k}(d)
\psi\left(\frac{x}
{d+\delta}\right)
\ll x^{11/24+7\varpi/12+\varepsilon}.
\]

By Lemma \ref{z3}, we have
\begin{align}\label{c20}
\begin{split}\sum_{D<d\leq 2D}&\tau_{k}(d)\psi\left(\frac{x}{d+\delta}\right)\\
&\ll
\left|\sum_{1\leq h\leq H}\frac1h\sum_{D<d\leq 2D}\tau_{k}(d) e\left(\frac{hx}{d+\delta}\right)\right|
\\&+ \left|\sum_{1\leq h\leq H}\frac1H\sum_{D<d\leq 2D}\tau_{k}(d)
e\left(\frac{hx}{d+\delta}\right)\right|
+D/H.
\end{split}\end{align}
Then we will focus on the estimate of
\begin{align*}
\sum_{1\leq h\leq H}\frac1h\sum_{D<d\leq 2D}\tau_{k}(d)e\left(\frac{hx}{d+\delta}\right).
\end{align*}
And we can deal with the sum
\[
\sum_{1\leq h\leq H}\frac1H\sum_{D<d\leq 2D}\tau_{k}(d)
e\left(\frac{hx}{d+\delta}\right)
\]
similarly.
By partial summation, we have
\begin{align*}
\begin{split}
&\sum_{D<d\leq 2D}\sum_{1\leq h\leq H}\frac{\tau_{k}(d)}{h}
e\left(\frac{hx}{d}\right)
e\left(\frac{-\delta hx}{d(d+\delta)}\right).
\\&
\ll \int_{D}^{2D}e\left(\frac{-\delta hx}{u(u+\delta)}\right)d\left(
\sum_{D<d\leq u}\sum_{1\leq h\leq H}\frac{\tau_{k}(d)}{h}
e\left(\frac{hx}{d}\right)\right)\\
&\ll
\max_{D\leq D_{1}\leq 2D}\left|\sum_{D<d\leq D_{1}}\sum_{1\leq h\leq H}\frac{\tau_{k}(d)}{h}
e\left(\frac{hx}{d}\right)\right|\\
&+\frac{Hx}{D^{2}}\max_{D< D_{1}\leq 2D}\left|\sum_{D<d\leq D_{1}}\sum_{1\leq h\leq H}\frac{\tau_{k}(d)}{h}
e\left(\frac{hx}{d}\right)\right|.
\end{split}
\end{align*}
We choose  $H=D^{2}/x^{1-\varrho}.$ Hence we have $H\geq 1.$
Then we will focus on the estimate of
\begin{align*}
S_{0}:=\sum_{D<d\leq D_{1}}\sum_{1\leq h\leq H}\frac{\tau_{k}(d)}{h}
e\left(\frac{hx}{d}\right)\end{align*}
By using the relation
$$\sum_{n_{1}n_{2}\cdots n_{k}=n}1=\tau_{k}(n),$$
and the dichotomy  method, we have
\[
\sum_{D<d\leq D_{1}}\sum_{1\leq h\leq H}\frac{\tau_{k}(d)}{h}
e\left(\frac{hx}{d}\right)\ll D^{\varepsilon}\sum_{1\leq h\leq H}\frac1h\sum_{d_{i}\sim D_{i},i=1,2,\cdots, k}e\left(\frac{hx}{d_{1}d_{2}\cdots d_{k}}\right),
\]
where
\begin{align}\label{k1}
d_{i}\leq d_{i+1}, \ D_{i}\leq D_{i+1}, \ \textup{for} \  1\leq i \leq k-1\end{align}
and
\begin{align}
\label{k2}
\prod_{i=1}^{k}D_{i}\sim D.\end{align}

Now we divide three cases to deal with this.

{\bf Case I}

Suppose that $D_{k}\geq D^{2/3}.$
Similar arguments for $T_{1}$ and $T_{2}$ in section 2, we have
$
S_{0}\ll x^{7/15 +\varepsilon},
$
where  we have chosen   $H=D^{2}/x^{1-\varrho}.$

{\bf Case II}

Suppose that  $D^{1/3}\leq D_{k}\leq D^{2/3}.$
By choosing $\mathcal{N}=x^{7/15}$ and  $(\kappa,\lambda)=(1/2,1/2)$ in Lemma \ref{LWY}, and restricted the range to $D^{1/3}\leq D_{k}\leq D^{3/8}$ and $D^{3/8}\leq D_{k}\leq D^{1/2},$  according to the symmetry, similar arguments as the argument of section 2 of $T_{3}$ and $T_{4}$, we have
\begin{align*}
S_{0}\ll
D^{\varepsilon}\sum_{1\leq h\leq H}\frac1h\sum_{d_{i}\sim D_{i},i=1,2,\cdots, k-1}\sum_{d_{k}\sim D_{k}}e\left(\frac{hx}{d_{1}d_{2}\cdots d_{k}}\right) \ll x^{7/15  +\varepsilon}.
\end{align*}

 {\bf Case III}

 Suppose that $D_{k}\leq D^{1/3}.$
 Then by (\ref{k1}) and (\ref{k2}), we have $D_{i}\leq D^{1/3},$ $i=1,2,\ldots,k.$ We   also suppose that $t$ is the least integer such that $D_{1}D_{2}\ldots D_{t}>D^{1/3}.$ Then we have
 \[
 D^{1/3}\leq (D_{1}D_{2}\ldots D_{t-1})D_{t}\leq D^{2/3}.
 \]
 Let $l_{1}=d_{1}d_{2}\ldots d_{t}$ and let $l_{2}=d_{t+1}d_{t+2}\ldots d_{k}.$
 Then we have
\[
S_{0}\ll
D^{\varepsilon}\sum_{1\leq h\leq H}\frac1h\sum_{l_{1}\sim L_{1} }\tau_{t}(l_{1})\sum_{l_{2}\sim L_{2} }\tau_{k-t}(l_{2})e\left(\frac{hx}{l_{1}l_{2} }\right),
\]
where $D^{1/3}\leq L_{1}\leq D^{2/3}$
 and $D^{1/3}\leq L_{2}\leq D^{2/3}.$
Then similar as the second case (consider $D^{1/3}\leq L_{1}\leq D^{3/8}$ and $D^{3/8}\leq L_{1}\leq D^{1/2}$ respectively), we have
$
S_{0} \ll x^{7/15 +\varepsilon}.
$

Then from the above three cases, we have
\[S_{0}:=\sum_{1\leq h\leq H}\frac1h\sum_{D<d\leq D_{1}}\tau_{k}(d)e\left(\frac{hx}{d}\right)\ll x^{7/15 +\varepsilon}.
\]
Then by (\ref{c100000})-(\ref{c20}), we have
\[S_{\tau_{k}}(x)= x\sum_{d=1}^{\infty}\frac{\tau_{k}(d)}{d(d+1)}
 +O\left(x^{7/15+\varrho+\varepsilon}
+x^{11/24+7\varpi/12+\varepsilon}+
x^{1/2-\varpi-\varrho+\varepsilon}
\right).\]
Choose
$\varrho=1/195$ and $\varpi=3/130.$
Recall   (\ref{c110}) and (\ref{c120}), then we can finally give Theorem \ref{th1}.

\section{Proof of Theorem \ref{th3}}
Let
\[
\mathcal{N}=x^{107/229}.
\]
We can write
\[
S_{\tau_{(2)}}(x):=S_{\tau_{(2)},1} +S_{\tau_{(2)},2} ,
\]
where
\begin{align}\label{c1111}
S_{\tau_{(2)},1}=\sum_{1\leq n\leq \mathcal{N}}
\tau_{(2)}
\left(\left[\frac{x}{n}\right]\right)
\end{align}
and
\begin{align}\label{c1211}
S_{\tau_{(2)},2}=\sum_{\mathcal{N}<n\leq x}\tau_{(2)}
\left(\left[\frac{x}{n}\right]\right).
\end{align}
Obviously, by   $\tau_{(2)}(n)\ll n^{\varepsilon},$ we have
\begin{align*}
S_{\tau_{(2)},1}=\sum_{1\leq n\leq \mathcal{N}}\tau_{(2)}\left(\left[\frac{x}{n}\right]\right)
&=\sum_{n\leq \mathcal{N}}(x/n)^{\varepsilon}
\\
&\ll \mathcal{N}^{1+\varepsilon}\\&\ll x^{ 107/229+\varepsilon}.
\end{align*}

As to $S_{\tau_{(2)},2},$  by    $\tau_{(2)}(n)\ll n^{\varepsilon},$
  we have
$$
\sum_{1\leq n\leq x} \tau_{(2)}(n)  \ll x^{ 1+\varepsilon}.
$$
Hence we can get
\begin{align}\label{c111}
\begin{split}
S_{\tau_{(2)},2}&=\sum_{\mathcal{N}<n\leq x}\tau_{(2)}\left(\left[\frac{x}{n}\right]\right)
\\&
=\sum_{d\leq x/\mathcal{N}}\tau_{(2)}(d)\sum_{x/(d+1)<n\leq x/d}1
\\&
=\sum_{d\leq x/\mathcal{N}}\tau_{(2)}(d)\left(\frac{x}{d}-\frac{x}{d+1}
-\psi\left(\frac{x}{d}\right)
+\psi\left(\frac{x}{d+1}
\right)\right)
\\&=x\sum_{d=1}^{\infty}
\frac{\tau_{(2)}(d)}{d(d+1)}
+O\left(
\mathcal{N}^{1+\varepsilon}
\right)
\\
&+O\left((\log x)\max_{\mathcal{N}<D\leq x^{1/2+\varpi}}\left|\sum_{D<d\leq 2D}\tau_{(2)}(d)
\psi\left(\frac{x}
{d+\delta}\right)\right|\right)
\\
&+O\left((\log x)\max_{x^{1/2+\varpi}<D\leq x/\mathcal{N}}\left|\sum_{D<d\leq 2D}\tau_{(2)}(d)
\psi\left(\frac{x}{d+\delta}\right)\right|\right),
\end{split}\end{align}
where  $\delta\in\{0,1\}.$
For $\mathcal{N}<D\leq x^{1/2+\varpi},$ by Lemma \ref{former11}, we have
\[
\sum_{D<d\leq 2D}\tau_{(2)}(d)
\psi\left(\frac{x}
{d+\delta}\right)
\ll x^{11/24+7\varpi/12+\varepsilon}.
\]
Then for $x^{1/2+\varpi}<D\leq x/\mathcal{N},$ we need to estimate
\[\sum_{D<d\leq 2D}\tau_{(2)}(d)\psi\left(\frac{x}{d+\delta}\right).
\]
By Lemma \ref{z3}, we have
\begin{align}\label{c211}
\begin{split}\sum_{D<d\leq 2D}&\tau_{(2)}(d)\psi\left(\frac{x}
{d+\delta}\right)\\
&\ll
\left|\sum_{1\leq h\leq H}\frac1h\sum_{D<d\leq 2D}\tau_{(2)}(d)e\left(\frac{hx}{d+\delta}\right)
\right|
\\&+ \left|\sum_{1\leq h\leq H}\frac1H\sum_{D<d\leq 2D}\tau_{(2)}(d)e\left(\frac{hx}{d+\delta}\right)
\right|
+D/H.
\end{split}\end{align}
Then we will focus on the estimate of
\begin{align*}
S_{\delta}:=\sum_{1\leq h\leq H}\frac1h\sum_{D<d\leq 2D}\tau_{(2)}(d)e\left(\frac{hx}{d+\delta}\right).
\end{align*}
And we can estimate the sum
\[
\sum_{1\leq h\leq H}\frac1H\sum_{D<d\leq 2D}\tau_{(2)}(d)e\left(\frac{hx}{d+\delta}\right)
\]
similarly.
By partial summation, we have
\begin{align}\label{n511}
\begin{split}
&\sum_{D<d\leq 2D}\sum_{1\leq h\leq H}\frac{\tau_{(2)}(d)}{h}
e\left(\frac{hx}{d}\right)
e\left(\frac{-\delta hx}{d(d+\delta)}\right)
\\&
\ll \int_{D}^{2D}e\left(\frac{-\delta hx}{u(u+\delta)}\right)d\left(
\sum_{D<d\leq u}\sum_{1\leq h\leq H}\frac{\tau_{(2)}(d)}{h}
e\left(\frac{hx}{d}\right)\right)\\
&\ll
\max_{D<D_{1}\leq 2D}\left|\sum_{D<d\leq D_{1}}\sum_{1\leq h\leq H}\frac{\tau_{(2)}(d)}{h}
e\left(\frac{hx}{d}\right)\right|\\
&+\frac{Hx}{D^{2}}\max_{D<D_{1}\leq 2D}\left|\sum_{D<d\leq D_{1}} \sum_{1\leq h\leq H}\frac{\tau_{(2)}(d)}{h}
e\left(\frac{hx}{d}\right)\right|.
\end{split}
\end{align}
We choose  $H=D^{2}/x^{1-\varrho}.$ Hence we have $H\geq 1.$
Hence, we only need to estimate the sum \begin{align*}
\sum_{1\leq h\leq H}\frac1h\sum_{D<d\leq D_{1}}\tau_{(2)}(d)
e\left(\frac{hx}{d}\right).
\end{align*}
By Lemma \ref{rs11}, we have
\begin{align*}
&\sum_{1\leq h\leq H}\frac1h\sum_{D<d\leq D_{1}}\tau_{(2)}(d)
e\left(\frac{hx}{d}\right)\\
&\ll
\sum_{1\leq h\leq H}\frac1h\sum_{D<n_{1}n_{2}\leq D_{1}}
\tau(n_{1})g(n_{2})
e\left(\frac{hx}{n_{1}n_{2}}\right)\\
&\ll
\sum_{1\leq h\leq H}\frac1h\sum_{D<m_{1}m_{2}n_{2}\leq D_{1}}
g(n_{2})
e\left(\frac{hx}{m_{1}m_{2}n_{2}}\right)
\\
&\ll
\sum_{n_{2}\leq D_{1}}g(n_{2})\sum_{1\leq h\leq H}\frac1h\sum_{D/n_{2}<m_{1}m_{2}\leq D_{1}/n_{2}}
e\left(\frac{hx}{m_{1}m_{2}n_{2}}\right)\\
& \ll
\sum_{n_{2}\leq D_{1}}g(n_{2})\sum_{1\leq h\leq H}\frac1h\sum_{m_{1}\leq (D/n_{2})^{1/2}}\sum_{D/m_{1}n_{2}<m_{2}\leq D_{1}/m_{1}n_{2}}
e\left(\frac{hx}{m_{1}m_{2}n_{2}}\right)
\\
&\ll
x^{\varepsilon}\sum_{n_{2}\leq D_{1}}g(n_{2})\left(
\left(\frac{xD^{3/2}}{n_{2}^{5/2}}\right)^{1/4}
+\frac{D^{3/4}}{n_{2}^{3/4}}
+\frac{D^{2}}{n_{2}x}
\right)\\
&\ll x^{\varepsilon}
\left(x^{1/4}D^{3/8}+D^{3/4}+D^{2}/x\right).
\end{align*}
This gives that
\[
S_{\tau_{(2)},2} \ll
x^{\varepsilon}\left(x^{11/24+7\varpi/12}
+x^{1/4+\varrho}D^{3/8}+x^{\varrho}D^{3/4}+
x^{\varrho}D^{2}/x+D/H\right).
\]
We choose $\varrho=8/458$ and $\varpi=7/458$
and  $\mathcal{N}=x^{107/229}.$
Then by (\ref{c111})-(\ref{c211}), we have
\[S_{\tau_{(2)}}(x)= x\sum_{d=1}^{\infty}\frac{\tau_{(2)}(d)}{d(d+1)}
 +O\left(x^{5/11+8\varrho/11+\varepsilon}
+x^{11/24+7\varpi/12+\varepsilon}+
x^{1/2-\varpi-\varrho+\varepsilon}
\right).\]
Recall   (\ref{c1111}) and (\ref{c1211}), then we can finally give Theorem \ref{th3}.
This completes the proof.

\bigskip
{\bf Acknowledgements} The author would like to thank the the referee's significant effort for giving some  detailed corrections and suggestions, which makes the description of the paper more rigorous and clear.

\address{Wei Zhang\\ School of Mathematics and Statistics\\
               Henan University\\
               Kaifeng  475004, Henan\\
               China}
\email{zhangweimath@126.com}

\end{document}